\documentclass[a4paper,12pt]{article}
\font\header=cmssdc10 at 20pt
\font\subheader=cmssdc10 at 15pt
\usepackage[ansinew]{inputenc}   
\usepackage{amsfonts}
\usepackage{indentfirst}
  \usepackage[normalem]{ulem}
  \usepackage{graphicx}
  \usepackage{amssymb}
 \usepackage{draftcopy}
  \usepackage{amsmath}
  \usepackage{rotating}
  \usepackage{color}
  \usepackage{setspace}
  \usepackage{fancybox}
 \usepackage{fancyvrb}
	\usepackage{fancyhdr}
	\usepackage[top=0.5cm,left=1cm,right=1cm,bottom=0.5cm]{geometry}
  \usepackage{wrapfig}
    \usepackage{float}
  \usepackage{type1cm}
\usepackage{eso-pic}
\usepackage{color}
 \usepackage{amsthm,amsmath,amsfonts}
\usepackage{fullpage,enumerate}
\usepackage{tikz}
\usepackage[colorlinks,citecolor=blue,urlcolor=blue]{hyperref}
\usepackage{hypernat}
\usepackage{graphicx}

\usepackage{textcomp}

\usepackage{bbm}

\newcommand{\la}{\lambda}

\def\aa{\alpha}
\def\bb{\beta}
\def\gg{\gamma}

\def\la{\lambda}

\def\1{{\bf 1}}
\def\rf{\rfloor}
\def\lf{\lfloor}
\def\R{\mathbb R}
\def\N{\mathbb N}
\def\Z{\mathbb Z}

\def\E{\mathcal E}
\def\T{\mathcal T}

\def\bx{\bar X}
\def\by{\bar Y}

\begin{document}

{\header Metastability of a random walk with catastrophes}

\vskip1cm

Luiz Renato Fontes and Rinaldo B. Schinazi

\bigskip

Instituto de Matem\'atica e Estatística,
Universidade de S\~ao Paulo,
Rua do Mat\~ao 1010,
05508-090 São Paulo SP 
Brasil; email: lrfontes@usp.br

 Department of Mathematics, University of Colorado, Colorado Springs, CO 80933-7150, USA; e-mail: 
Rinaldo.Schinazi@uccs.edu

\vskip1cm

{\bf Abstract.} We consider a random walk with catastrophes which was introduced to model population biology. 
It is known that this Markov chain gets eventually absorbed at $0$ for all parameter values. Recently, it has been shown that this  chain exhibits a metastable behavior 
in the sense that it can persist for a very long time before getting absorbed. In this paper we study this metastable phase by making the parameters converge to extreme values.
We obtain four different limits that we believe shed light on the metastable phase.

\vskip1cm


\vskip1cm

{\header 1 A random walk with catastrophes}

\vskip1cm

Let $(X_n)_{n\geq 0}$ be the following discrete-time Markov chain. For every $n\geq 1,$ $X_n$ is a non-negative integer.
The model has two parameters, $p\in (0,1)$ and $c\in (0,1).$ For $n\geq 0,$ let $X_n=k\geq 0.$ If $X_n=0$ then $X_{n+1}=0$. That is,
$0$ is an absorbing state. If $k\geq 1$,  there are two possibilities:
\begin{itemize}
\item With probability $p$ there is a birth. Then, $X_{n+1}=k+1$.
\item With probability $1-p$ there is a catastrophe. Then, $X_{n+1}=k-B_n,$ where $B_n$ is a binomial random variable with
parameters $k$ and $c.$ The random variables $B_n$ are sampled independently of each other and of everything else.
\end{itemize}

Ben-Ari et al. (2019) have shown that this Markov chain is eventually absorbed at $0$ for all values of $p$ and $c$ in $(0,1)$.
They have also shown through simulations and first moment computations that the time to absorption can be unusually long, in particular if $p$ is close to 1
or $c$ is close to $0$. Before absorption the chain fluctuates in a narrow band around 
$$n^*=\frac{p}{(1-p)c}.$$
That is, the chain seems to have reached some equilibrium but this equilibrium turns out to be unstable, see Figure 1. This is why we think of the chain as exhibiting metastable behavior. By looking at the limits 
of the process  for $c\to 0$ and/or $p\to 1$ we will get more insight into this metastable phase. 
We find four different limit processes explicitly. 
In the expression of $n^*$, $(1-p)$ and $c$ play the same role. However, we will show that as $c$ approaches $0$ and $1-p$ approaches $0$ the limiting processes are very different.

This model goes back to at least Neuts (1994), see Section 2 there. The catastrophe distribution need not be binomial, see Brockwell (1986). For a more recent survey of these models, see Artalejo et al. (2007).

\begin{figure}[ht]
\centering
	\includegraphics[width=10cm]{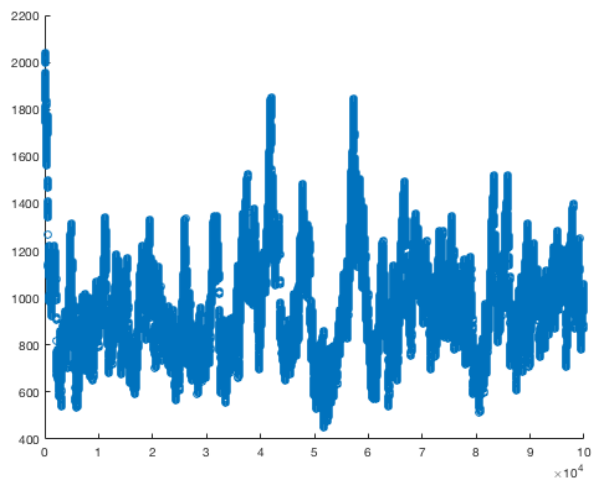}
	\caption{We ran this simulation for $10^5$ steps, starting at $X_0=2000$, with $p=0.99$ and $c=0.1$. We see that the process drops very fast to $n^*=990$ and then fluctuates around this value.}
	\end{figure}
	
	\vskip1cm

{\header 2 The limit as c approaches 0 and p is fixed}

\vskip1cm

Let $(Y_n)_{n\geq 0}$ be the following discrete-time Markov chain. Assume $Y_n=k\geq 1$,

\begin{itemize}
\item With probability $p$, $Y_{n+1}=k+1$.
\item With probability $1-p$, $Y_{n+1}=k-P_n,$ where $P_n$ is a Poisson random variable with mean $\frac{p}{1-p}$.
The random variables $P_n$ are sampled independently of each other and of everything else.
\end{itemize}

Let  $(X_n(L))_{n\geq 0}$ be the Markov chain defined in the previous section with $c=1/L$ and $p$ in $(0,1)$.
With  $c=1/L$ we get 
$n^*(L)=L\frac{p}{1-p}.$

\medskip

{\bf Proposition 1. } {\sl Let $X_0(L)=Y_0=n^*(L)$ then for any fixed $T$ and $p$,
$$\lim_{L\to+\infty}P(\exists n\leq T: X_n(L)\not=Y_n)=0.$$}

\medskip

In words, the process $(X_n(L))_{n\leq T}$ approaches $(Y_n)_{n\leq T}$ as $L$ goes to infinity.

\medskip

{\it Proof of Proposition 1}

We introduce an auxiliary process $(U_n)$. If $U_n=k\geq 1$, then
\begin{itemize}
\item With probability $p$, $U_{n+1}=k+1$.
\item With probability $1-p$, $U_{n+1}=k-P'_n,$ where $P'_n$ is a Poisson random variable with mean $kc$.
The random variables $P'_n$ are sampled independently of each other and of everything else.
\end{itemize}

We proceed in two steps. First we show that
$$\lim_{L\to+\infty}P(\exists n\leq T: X_n(L)\not=U_n)=0,$$
using that 
$$P(bin(x,c)\not =Poiss(xc))\leq \frac{1}{2}xc^2.$$
Then, we show that
$$\lim_{L\to+\infty}P(\exists n\leq T: U_n\not=Y_n)=0,$$
using that
$$P(Poiss(\lambda)\not =Poiss(\mu))\leq |\lambda-\mu|.$$

\medskip

First step.
Note that $U_n\leq n+n^*\leq T+n^*$ for $n\leq T$. Therefore, each time the process drops it drops from a state smaller than $T+n^*$.
There are at most $T$ drops up to time $T$. Hence, the sum of all drops up to time $T$ is stochastically less than a sum of $T$ i.i.d. Poisson variables with mean $(T+n^*)c$. The latter sum being itself a Poisson random variable with mean $T(T+n^*)c$.

Let 
$$A=\{Poiss\left(T(T+n^*)c\right)<T(T+n^*)c+M\},$$
where $M$ will be large.
Let $\tau=\min\{n\geq 1: X_n(L)\not=U_n\}.$ We have,
$$P(\tau=n)\leq P(X_{n-1}(L)=U_{n-1}; X_n(L)\not=U_n;A)+P(A^c).$$
On $A$,   the process $(U_{n})_{n\leq T}$ is a positive integer in the interval 
$$I=[n^*-T(T+n^*)c-M, T+n^*].$$
Hence,
$$P(X_{n-1}(L)=U_{n-1}; X_n(L)\not=U_n;A)\leq \sum_{x\in I}P(bin(x,c)\not=Poiss(xc)).$$
Therefore,
\begin{align*}
P(X_{n-1}(L)=U_{n-1}; X_n(L)\not=U_n;A)&\leq  \sum_{x\in I}\frac{1}{2}xc^2\\
&=\frac{1}{2}c^2\frac{1}{2}\left(2n^*+T-T(T+n^*)c-M\right)\left(T+T(T+n^*)c+M+1\right).
\end{align*}
Using that $cn^*$ is a constant and that $c$ goes to $0$ as $L$ goes to infinity,  
$$\lim_{L\to\infty}P(X_{n-1}(L)=U_{n-1}; X_n(L)\not=U_n;A)=0.$$
Hence,
$$\limsup_{L\to\infty} P(\tau=n)\leq P(Poiss(\lambda T)>\lambda T+M),$$
where $\lambda=cn^*$.
Letting now $M$ go to infinity we get $\lim_{L\to\infty} P(\tau=n)=0$ for every $n\leq T$. This completes the first step.

\medskip

Second step. Using the notation introduced in the first step,
$$P(U_{n-1}=Y_{n-1}; U_n\not=Y_n;A)\leq \sum_{x\in I}P(Poiss(xc)\not=Poiss(\lambda)).$$
Therefore,
$$P(U_{n-1}=Y_{n-1}; U_n\not=Y_n;A)\leq  \sum_{x\in I}|cx-\lambda|.$$
For $x$ in $I$,
$$-T(Tc+\lambda)c-Mc\leq cx-\lambda \leq Tc.$$
Hence, for $x$ in $I$ and $M$ large enough,
$$|cx-\lambda|\leq (M+T\lambda)c.$$
Therefore,
$$P(U_{n-1}=Y_{n-1}; U_n\not=Y_n;A)\leq \left(T(T+n^*)c+M+T+1\right)(M+T\lambda)c.$$
As $L$ goes to infinity the r.h.s. goes to 0 and we can conclude as in Step 1. This completes the proof of Proposition 1.

\vskip1cm

{\header 3 The limit as p approaches 1 and c is fixed}

\vskip1cm

We now make $p=1-\frac1L$ and let $c\in(0,1)$ be fixed. Let $y\geq0$ and set $X_0=\lf yL\rf$.
In this context, it is convenient to describe $(X_n)$ as follows. 
It alternates between the following two modes. 
\begin{itemize}
	\item  Forward mode (F): The chain jumps to the right a geometric number of steps. Each jump takes a unit time. The mean of the geometric random variable is $L-1$.     
	\item  Backward mode (B):  The chain jumps to the left in one unit time from the position it has gotten to, say $z_L$. The size of the jump is distributed 
    as a Binomial random variable with mean $cz_L$ and variance $c(1-c)z_L$.
\end{itemize}
The chain starts with Mode (F) with probability $p$ or with Mode (B) with probability $1-p$. Then it alternates deterministically between the two modes.

\medskip

We now introduce what will turn out to be the limiting process of the rescaled process $(\frac{1}{L}X_{tL})_{t\geq 0}$. 

Let $\E_0,\E_1,\ldots$ be independent mean 1 exponential 
random variables. Let $S_0=0$, and $S_n=\sum_{i=1}^n\E_i$, $n\geq1$. Define recursively $Y_0=y$ and, for $n\geq1$, 
\begin{align*}
Y_{S_n-}&=Y_{S_{n-1}}+\E_n \\
Y_{S_n}&=(1-c)Y_{S_n-}
\end{align*}
 
For $t$ in $(S_{n-1}, S_n)$, 
let $Y_t$ be the  linear interpolation of $Y_{S_{n-1}}$ and $Y_{S_n}$.

In words, starting from $y$, $(Y_t)_{t\geq0}$ first moves to the right at speed 1 for a mean 1 exponential
random time, after which it finds itself at $y+\E$, and then it jumps instantaneously to the left by 
$c(y+\E)$ units. Forward and backward jumps keep alternating in a deterministic way.

\medskip

{\bf Proposition 2. } {\sl Let $p=1-\frac1L$ and let $c$ fixed in $(0,1)$ be fixed. Let $y\geq 0$ and set $X_0=\lf yL\rf$. Then,
$\left(\frac1LX_{tL}\right)_{t\geq0}$ (with the proper interpolated definition of $X_t$ for $t$ outside $\frac1L\N$) converges weakly as $L\to\infty$ to $(Y_t)_{t\geq0}$ where $Y_0=y$.}

\medskip

{\it Proof of Proposition 2}

As $L$ goes to infinity, $p$ goes to 1. So the first mode taken by the chain after time $0$  is (F). Therefore, let $n\in \frac1L\N$ then
$X_{nL}=\lf yL\rf+nL$ for $nL\leq G_1$ where $G_1$ is a mean $L-1$ geometric random variable. Since
$G_1/L$ converges weakly to $\E_1$ (a mean 1 exponential random variable) we get the following weak convergence,
$$\lim_{L\to+\infty} \frac{1}{L} X_{nL}=y+n\mbox{ for }n<\E_1.$$
This shows the existence of $\lim_{L\to+\infty} \frac{1}{L} X_{tL}$ for $t$ in  $\frac1L\N$. For $t$ outside this set we extend the definition of  $X_{tL}$ by linear
interpolation. This gives 
$$\lim_{L\to+\infty} \frac{1}{L} X_{tL}=y+t \mbox{ for all }t<\E_1^{-}.$$

After the first mode (F) we switch to mode (B).
Let $LT_L=G_1$,  note that $T_L$ converges weakly as $L$ goes to infinity is $\E_1$.
Let $z_L=X_{LT_L}$ then
 $$X_{LT_L+1}=z_L-bin(z_L,c).$$
 Since $z_L$ goes to infinity with $L$, by the Law of Large Numbers and the weak convergence of $G_1/L$ we get the following weak convergence,
 $$\lim_{L\to+\infty}\frac{1}{L}X_{LT_L+1}=(1-c)(\E_1+y).$$
 At this point we have shown that  $(\frac{1}{L} X_{tL})$ converges to $(Y_t)$ for the first (F) and (B) modes.
Using this method we can continue computing limits for the successive (F) and (B) modes. Since the latter process is non explosive, convergence in the usual Skorohod space of c\`adl\`ag trajectories readily follows.

\medskip

{\bf Proposition 3. }{\sl The process $(Y_t)_{t\geq 0}$ is ergodic. That is, it converges weakly to an invariant measure. 
Moreover, the distribution of the invariant measure is the same as the distribution of $
	\sum_{n\geq0}(1-c)^n\E_n,
	$
	where $\E_0,\E_1,\ldots$ are mean 1 i.i.d. exponential random variables}
	
	\medskip
	
	{\it Proof of Proposition 3}

    The infinitesimal generator of $Y$ for $f\in C^1_c$, the continuously differentiable real functions on $\R^+$ with compact support, is given by
    \begin{equation}\label{gen}
    	\Omega f(x) = f'(x) + f((1-c)x) - f(x).
    \end{equation}
    To justify this, write
    $$E_x[f(Y_t)] = f(x+t)e^{-t}+\int_0^t ds\,e^{-s} f((1-c)(x+s)+t-s)+o(t);$$
    thus,
    $$E_x[f(Y_t)]-f(x) = [f(x+t)-f(x)]e^{-t}+\int_0^t ds\,e^{-s} [f((1-c)(x+s)+t-s)-f(x)]+o(t),$$
    and~(\ref{gen}) follows by dividing by $t$ and taking the limit as $t\to0$.
    
    \medskip

    In order to find an invariant distribution, let us suppose one such measure admits a continuous density $\psi$, which thus must satisfy
    \begin{align*}
    	\int_0^\infty\Omega f(x) \psi(x) dx &= \int_0^\infty f'(x) \psi(x) dx - 
    	\int_0^\infty dx\,\psi(x) \int_{(1-c)x}^x f'(y)dy\\
    	&= \int_0^\infty f'(x) \left\{\psi(x) - \int_{x}^{ax} \psi(y)dy\right\}=0
    \end{align*}
   for all $f\in C^1_c$, where $a=(1-c)^{-1}$. It follows that
   \begin{align*}
	  \psi(x) = \int_{x}^{ax} \psi(y)dy
    \end{align*}
    for all $x>0$. By taking Laplace transforms, we readily find that 
    $\varphi(\theta)=\int_0^\infty e^{-\theta x} \psi(x) dx$, $\theta>0$, must satisfy
    $$\varphi(\theta)=\frac1{1+\theta}\varphi((1-c)\theta),\,\theta>0.$$
    Iterating and taking the appropriate limit, we find that
    $$\varphi(\theta)=\prod_{n=0}^\infty\frac1{1+(1-c)^n\theta},\,\theta>0,$$
    and the claimed form of the invariant measure is established in this case. 
    
    To verify uniqueness, we resort to Meyn and Tweedie (1993). 
    \begin{enumerate}
    	\item (Lebesgue)-irreducibility: It is enough to check the condition (midway at page 490) for $B$ a finite nonempty open interval $(a,b)$ with $a>0$. The condition is clear for $x<b$; if $x\geq b$, then it is enough 
    	to establish that $P_x(\tau_{(0,a)}<\infty)>0$, but this follows from the fact that after $n$ jumps, our process is found at 
    	\begin{equation}\label{wn}
    		W_n:=(1-c)^nx+\sum_{i=1}^n(1-c)^{n-i+1}\,\E_i,
    	\end{equation}
    with $\E_1,\E_2,\ldots$ as above. It is enough now to have $n$ large enough to make the first term less than $a/2$ and then $\E_1,\ldots,\E_n$ small so as to make the sum in the second term less than $a/2$, and event of positive probability.
    	\item Non-evanescence: If $Y_t\to\infty$ as $t\to\infty$, then $W_n\to\infty$ as $n\to\infty$, with $W_n$ as
    	in~(\ref{wn}). But $W_n$ is stochastically bounded by $Z:=x+\sum_{n\geq1}(1-c)^n\E_n$ uniformly in $n$. Since $Z$ is a proper random variable, it follows that $P_x(\lim_{t\to\infty}Y_t=\infty)=0$ for all $x$.
    	\item T-process property: We resort to Theorem 4.1 of \textlangle M-T\textrangle. We already have non-evanescence,
    	so we want to argue that $[0,L]$ is {\it petite} for every $L>0$; we want to exhibit a probability measure $a=a_L$ and a nontrivial measure $\nu=\nu_L$ on $\R^+$ such that $\int_0^\infty da(t) P_x(Y_t\in\cdot)\geq\nu(\cdot)$ for all $x\in[0,L]$.
    	We choose $a(t)=e^{-t+L}{\mathbbm 1}\{t>L\}$ and thus
    	\begin{align*}
    	&\int_0^\infty da(t) P_x(Y_t\in\cdot)=\int_L^\infty dt\,e^{-t+L} P_x(Y_t\in\cdot)\\
    	\geq& \int_L^\infty dt\,e^{-t+L} P_x(\E_1>L-x,\,Y_t\in\cdot)=\int_L^\infty dt\,e^{-t+L} e^{-L+x} P_L(\,Y_{t-L+x}\in\cdot)\\
    	=&\int_x^\infty dt\,e^{-s-L+2x} P_L(\,Y_{s}\in\cdot)\geq e^{-L}\int_L^\infty ds\,e^{-s}\,P_L(Y_s\in\cdot)=:\nu(\cdot)
    	\end{align*}
    	for all $x\in[0,L]$, and we have found our measure $\nu$.
    	\item To conclude, we apply Theorem 3.2 of Meyn and Tweedie (1993) to get that $X$ is Harris recurrent. This is sufficient for uniqueness of the invariant distribution, as pointed out in Meyn and Tweedie (1993)--- (see last but one sentence of the second paragraph in page 491).
   \end{enumerate}

\medskip

{\it Remark.} One amusing point related to the above proposition is as follows. The discrete time processes $M_n:=Y_{S_n-}$ and 
$m_n:=Y_{S_n}$, $n\geq1$, represent local maxima and minima of the trajectory of $Y$. One would then perhaps be led to guess that the invariant distribution of $Y$ should (strictly) dominate the invariant distribution of $m_n$, and be dominated by the invariant distribution of $M_n$. However, it {\it equals} the latter distribution (as one may easily check by computing the invariant distributions of $m_n$ and $M_n$). The apparent contradiction is dispelled by the realization that (looking at the invariant distribution of $Y$ as the limiting distribution of $Y_t$ as $t\to\infty$) the interval $(m_n, M_{n+1})$ containing $t$ is larger than typical (this is of course an instance of the {\it inspection paradox}), and in this case it is asymptotically 'twice' the size of a typical interval. Indeed, one can argue along this line to show that $\lim_{t\to\infty}Y_t=m_\infty+\E=M_\infty$ in distribution, where $m_\infty$ and $M_\infty$ are the invariant distributions of $m_n$ and $M_n$ respectively, with $m_\infty$ and $\E$ independent; actually, this provides an alternative proof of Proposition 3.
    
%
%

\vskip1cm

{\header 4 The limit as p approaches 1 and c approaches 0}

\vskip1cm

We now make $p=1-\frac1{L^\aa}$ and let $c=\frac1{L^\bb}$, $\aa\in(0,1)$, $\bb=1-\aa$, 
and let $X_0=R_L+\lf yL^\aa\rf$, with $R_L=rL+o(L)$,
$r\in\R^+$, 
$y\in\R$ fixed, and set $\bx_t=X_t-R_L$.

The limiting process is defined by (using the notation of Section 3) $\by_0=y$ and, for $n\geq1$, $\by_{S_n-}=\by_{S_{n-1}}+\E_n$, $\by_{S_n}=\by_{S_n-}-r$, and 
linear interpolation on $(S_{n-1}, S_n)$. 

\medskip

{\bf Proposition 4. }{\sl The rescaled and centered process 
$\left(\frac1{L^\aa}\bx_{tL^\aa}\right)_{t\geq0}$ converges weakly as $L\to\infty$ to $(\by_t)_{t\geq0}$}.

\medskip

With our choice of parameters the metastable equilibrium $n^*$ is of order $L$. The initial state is of order $rL$. Maybe surprisingly the limiting process drifts linearly with a speed $1-r$.

\medskip

{\it Proof of Proposition 4}

We follow the analysis done in Proposition 2. Again we start as $L$ goes to infinity with a forward mode. 
Let $T'$ be such that $T'L^{\alpha}=G'$ where $G'$ is a geometric random variable with mean $L^\alpha-1$. Then,
$X_{T'L^{\alpha}}=X_0+G'$.
Since $X_0= R_L+\lf yL^\aa\rf$,
$$X_{T'L^{\alpha}}=R_L+\lf yL^\aa\rf+G'.$$
Therefore, $\bx_{T'L^{\alpha}}=\lf yL^\aa\rf+G',$
and
$$\lim_{L\to\infty}\frac1{L^\aa}\bx_{T'L^\aa}=y+\E_1,$$
where $\E_1$ is a mean 1 exponential random variable.

At time  $T'L^{\alpha}+1$ we switch to the backward mode with a single jump. Let $z'_L=X_{T'L^{\alpha}}$. Then,
$$X_{T'L^{\alpha}+1}=z'_L-bin(z'_L,\frac{1}{L^\beta}).$$
Hence,
$$\frac1{L^\aa}\bx_{T'L^\aa+1}=\frac 1{L^\aa}\left (\lf yL^\aa\rf+G'-bin(z'_L,\frac{1}{L^\beta})\right).$$
In order to prove that the limit of the l.h.s. is $y+\E_1-r$ we need to show that
$$\lim_{L\to +\infty}\frac 1{L^\aa}bin(z'_L,\frac{1}{L^\beta})=r.$$
We do this next. Note that $z'_L=rL+o(L)$ and let
$$
J_L:=\frac1{L^\aa}\sum_{i=1}^{rL+o(L)}\xi_i,
$$
where $\xi_1,\xi_2,\ldots$ are iid Bernoullis with parameter $1/L^\bb$. Taking the Laplace transform, we find
$$
E\left(e^{\theta J_L}\right)=\left( 1+\frac1{L^\bb}\left( e^{\theta/L^\aa}-1\right)\right)^{rL+o(L)}
=\left( 1+\frac\theta{L}\,\frac{e^{\theta/L^\aa}-1}{\theta/L^\aa}\right)^{rL+o(L)}\to e^{r\theta},
$$
for all $\theta$, and this shows that $J_L\to r$ in probability as $L\to\infty$.

Using this method we can continue computing limits for successive forward and backward modes.

\vskip1cm

{\header 5 The limit as p and c approach 0, c faster}

\vskip1cm

We now take $p=\frac1{L^\gg}$, $c=\frac1{L^{1+\gg}}$, $\gg>0$. Let us make $X_0=R_L+k$, with $R_L=rL+o(L)$,
$r\in\R^+$, $k\in\Z$ fixed, and set $\bx_t=X_t-R_L$. Notice that $\bx_0=k$.

\medskip

{\bf Proposition 5. }{\sl The rescaled and centered process $(\bx(tL^\gg))_{t\geq0}$ converges weakly as $L\to\infty$ to a continuous time simple random walk on $\Z$ with jump rate $1+r$, jumping to the right with probability $\frac1{1+r}$.}

\medskip

{\it Proof of Proposition 5}

%
%

Let us describe the jump times and sizes of $X_t$ starting at a location $Q_L=rL+o(L)$ as follows.
Let $Z^1_1,Z^1_2,\ldots$ be independent Bernoulli random variables with success parameter $1/L^\gg$, and, independently, let $Z^2_1,Z^2_2,\ldots$ be independent binomial random variables with $Q_L$ trials and success parameter $1/L^{1+\gg}$. 
Now set $T^i_L=\inf\{j\geq1:\,Z^i_j>0\}$, $i=1,2$. Notice that $T_L^1$ and $T_L^2$ are independent geometric random variables. Moreover, for $i=1,2$, as $L\to\infty$ we have the following convergence in distribution, $\frac1{L^\gg}T^i_L\to\T^i$ where $\T^i$ is a rate $\la_i$  exponential random variable and $\la_1=1$, $\la_2=r$. 

Hence, the time of the first jump is $T_L:=T_L^1\wedge T_L^2$. Note that  $\frac1{L^\gg}T_L$ converges in distribution to an exponential random variable with rate $\la_1+\la_2=1+r$.

We now turn to the the jump length.  If $T_L^1\leq T_L^2$ then the chain jumps one unit to the right.  As $L\to\infty$ this has probability 
$\frac1{1+r}$.  If  $T_L^1>T_L^2$ then the jump is equal $-Z^2_{T_L}$. Note that this is a strictly negative integer. We claim that
$$P(Z^2_{T_L}\geq2|T_L^1>T_L^2)\to 0$$
as $L\to\infty$. It is enough to show that $P(Z^2_{T^2_L}\geq2)\to 0$ as $L\to\infty$. The latter probability equals
$$\frac{P(Z^2_1\geq2)}{P(Z^2_1\geq1)}.$$
The denominator equals $1-\left(1-\frac1{L^{1+\gg}}\right)^{Q_L}\geq 1-e^{-\frac r{2L^\gg}}\geq\frac{r}{3L^\gg}$ for $L$ large enough; and the numerator
is bounded above by $(Q_L/L^{1+\gg})^2\leq 2r^2/L^{2\gg}$ for $L$ large enough. The claim is established. This shows that in the limit when the process  $(\bx(tL^\gg))_{t\geq0}$ jumps to the left it jumps exactly one unit. The proof of Proposition 5 is complete.

\medskip

{\it Remark.} A note about the distinction between $R_L$ and $Q_L$. The former quantity is part of the position of $X$ at time $0$, while the latter is meant for a generic position of the process after a fixed number of steps (independent of $L$) --- the $o(L)$ of $R_L$ is fixed, and might have been written as $o_0(L)$, while that of $Q_L$ varies from step to step.

%
%
%
%
%
%
%
%
%
\vskip1cm

{\subheader Acknowledgements}

LRF is supported in part by CNPq grant 311257/2014-3 and FAPESP grant 2017/10555-0. 

\vskip1cm

{\header References}

\vskip1cm

J.~R.~Artalejo, A.~Economou and M.~J.~Lopez-Herrero (2007). Evaluating growth measures in populations subject to
binomial and geometric catastrophes. Math. Biosci. Eng. 4, 573-594.

I. Ben-Ari, A. Roitershtein, R.B.Schinazi (2019) A random walk with catastrophes. Electron. J. Probab. 24, 1-21.

P.~J.~Brockwell (1986).
The extinction time of a general birth and death process with catastrophes.
J. Appl. Probab. 23, 851-858.

S.P. Meyn and R.L. Tweedie (1993) Stability of Markovian processes II: continuous time processes and sampled chains. Adv. Appl. Prob. 25, 487-517.

M.~F.~Neuts (1994). An interesting random walk on the non-negative integers.
J. Appl. Probab. 31, 48-58.

\end{document}